

\input epsf.tex

\font\tenfrak=eufm10
\font\sevenfrak=eufm7
\font\fivefrak=eufm5
\newfam\frakfam
\textfont\frakfam=\tenfrak
\scriptfont\frakfam=\sevenfrak
\scriptscriptfont\frakfam=\fivefrak

\font\tenbb=msbm10
\font\sevenbb=msbm7
\font\fivebb=msbm5
\newfam\bbfam
\textfont\bbfam=\tenbb
\scriptfont\bbfam=\sevenbb
\scriptscriptfont\bbfam=\fivebb

\font\tensym=msam10
\font\sevensym=msam7
\font\fivesym=msam5
\newfam\symfam
\textfont\symfam=\tensym
\scriptfont\symfam=\sevensym
\scriptscriptfont\symfam=\fivesym

\font\tenlas=lasy10
\font\sevenlas=lasy7
\font\fivelas=lasy5
\newfam\lasfam
\textfont\lasfam=\tenlas
\scriptfont\lasfam=\sevenlas
\scriptscriptfont\lasfam=\fivelas

\def\bbR{{\Bbb R}}
\def\bbZ{{\Bbb Z}}

\def\bbC{{\Bbb C}}
\def\bbQ{{\Bbb Q}}
\def\bbP{{\Bbb P}}

\def\bbF{{\Bbb F}}

\def\calD{{\Cal D}}

\def\calM{{\Cal M}}
\def\calN{{\Cal N}}

\def\calP{{\Cal P}}

\def\calA{{\Cal A}}

\def\calT{{\Cal T}}

\def\Pic{\roman{Pic}}
\def\Ima{\roman{Im}}

\def\med{\medskip\noindent}

\def\Aut{\roman{Aut}}
\def\Bir{\roman{Bir}}
\def\Kum{\roman{Kum}}

\def\Lee{\roman{Leech}}

\def\med{\bigskip\noindent}

\input amstex

\documentstyle{amsppt}
\NoBlackBoxes

\topmatter
\title
Birational automorphisms of quartic  Hessian surfaces
\endtitle

\author
Igor V. Dolgachev and Jonghae Keum
\endauthor

\abstract
We find generators of the group of birational automorphisms of the
 Hessian surface of a general cubic surface. Its nonsingular minimal
model is a K3 surface with the Picard lattice of rank 16. The latter  embeds
naturally in the even unimodular lattice
$II^{1,25}$ of rank 26 and signature $(1,25)$ as the orthogonal complement of a root sublattice of rank 10. Our generators are
related to reflections with respect to some Leech roots.  A similar
observation was made first in the case of quartic Kummer surfaces in the
work of S. Kond$\bar {\roman o}$. We shall explain   how
our generators are related to the generators of the group of birational
automorphisms of a general quartic Kummer surface which
is birationally isomorphic to a special Hessian surface.

\endabstract
\thanks Research of the first author was supported in part by  NSF
grant DMS 9970460.
Research of the second author was supported by
 Korea Research Foundation Grant
(KRF-2000-D00014).
\endthanks
\endtopmatter
\document

\head {\bf  0. Introduction}\endhead

Let $S:F(x_0,x_1,x_2,x_3) = 0$ be a nonsingular cubic surface in $\bbP^3$
over $\bbC$. Its Hessian surface is a quartic surface defined by the
determinant of the matrix of second order partial derivatives of the
polynomial
$F$. When
$F$ is general enough, the quartic $H$ is irreducible and has 10
nodes. It contains also 10 lines which are the intersection lines
of five planes in general linear position. The union of these five planes
is classically known as the Sylvester pentahedron of $S$. The equation of
$S$ can be written as the sum of cubes of some linear forms defining the
five planes. A nonsingular model of $H$ is a K3-surface $\tilde H$. Its
Picard number $\rho$ satisfies the inequality $\rho \ge 16$. In this paper
we give an explicit description of the group $\Bir(H)$ of birational
isomorphisms of $H$ when $S$ is general enough so that $\rho = 16$.
 Although $H$, in general, does not have any non-trivial
automorphisms (because $S$ does not), the group $\Bir(H)\cong
\Aut(\tilde H)$ is infinite.  We show that it is generated by the
automorphisms defined by projections from the nodes of $H$, a
birational involution which interchanges the nodes and the lines,
and the inversion automorphisms of some elliptic pencils on
$\tilde H$. This can be compared with the known structure of the
group of automorphisms of the Kummer surface associated to the
Jacobian abelian surface of a general curve of genus 2 (see
[Ke2], [Ko]).  The latter surface is birationally isomorphic to
the Hessian $H$ of a cubic surface ([Hu1]) but the Picard number
of $\tilde H$ is equal to 17 instead of 16. We use the method for
computing Bir$(H)$ employed by S. Kond$\bar {\roman o}$ in [Ko].
 We show that the Picard lattice $S_H$ of the K3-surface $\tilde
H$  can be primitively embedded into  the unimodular lattice $L =
\Lambda\perp U$, where $\Lambda$ is the Leech lattice and $U$ is
a hyperbolic plane. The orthogonal complement of $S_H$ in $L$ is
a primitive lattice of rank 10 which contains a root lattice of
type $A_5\perp A_1^{\perp 5}$. The fundamental domain $D$  of the reflection
group of $L$ cuts out a finite polyhedron $D'$ in the fundamental
domain $D_H$ for the (-2)-reflection group of $S_H$ in the
connected component of the set $\{x\in S_H\otimes \bbR: (x,x) >
0\}$ containing an ample divisor class. We determine the
hyperplanes $h$ which bound $D'$ and match them with
automorphisms $h\to \sigma_h$ of $\tilde H$ such that $\sigma_h$
sends one of the half-spaces defined by $h$ to the other
half-space defined by $h$ or to one of the two half-spaces corresponding to 
 $\sigma_h^{-1}$.
 This allows one to prove that the automorphisms
$\sigma_h$ generate a group of symmetries of $D_H$, having $D'$
as its fundamental domain. By Torelli Theorem for K3-surfaces this
implies that the automorphisms $\sigma_h$ together with some
symmetries of $D'$ generate the group Aut$(\tilde H)$. The reason
why the beautiful combinatorics of the Leech lattice plays a role
in the description of the automorphisms of $H$ is still unclear to
us. We hope that the classification of all K3-surfaces whose
Picard lattice is isomorphic to the orthogonal complement of a
root sublattice of $L$ will shed more light to this question.

\head  {\bf 1. Some classical facts}\endhead

\bigskip\noindent
Here we summarize without proofs some of the known properties of the
quartic Hessian surface
$H$ of a cubic surface $F$ in $\bbP^3$. We refer for the proof to the
classical literature (for example, [Ba,Ca]).
\med
\item{(1)} $H$ is the locus of points $x\in \bbP^3$ such that the
polar quadric $P_x(F)$ of $F$ is singular.
\item{(2)} $H$ is the locus of points $x\in \bbP^3$  such that there
exists a polar quadric $P_y(F)$ of $F$ such that $x\in
\roman{Sing}(P_y(F))$.
\item{(3)} If $F$ is nonsingular and general enough, then $H$ has 10
nodes corresponding
to polar quadrics of corank 2 and has 10 lines corresponding to their
singular lines.
\item{(4)} The  lines (resp. singular points) of $H$ are the
intersection lines of 10 pairs (resp. 10 triples) of hyperplanes
$\pi_i, i = 1,\ldots,5$, any four of them being linearly
independent. The union of the planes $\pi_i$ is called the Sylvester
pentahedron of the cubic surface.
\item{(5)} If $l_i = 0$ are the equations of the hyperplanes
$\pi_i$, the equation of $F$ can be written in the Sylvester form
$a_1l_1^3+\ldots+a_5l_5^3 = 0$. The equation of $H$ can be written in the
form
${1\over a_1l_1}+\ldots+{1\over a_5l_5} = 0$.
\item{(6)} $H\cap F$ is the parabolic curve of $F$: the set of points
$x\in F$ such that the tangent hyperplane $T_x(F)$  intersects $F$
along a cubic curve which has a cuspidal points at $x$.
\item{(7)} If $x\in H$ and $y\in \roman{Sing}(P_x(F))$ then the
correspondence $x\to y$ is a birational involution $\tau$ of $H$.
On a minimal nonsingular model $\tilde H$ it interchanges the
exceptional curve blown up from the node $P_{ijk} =
\pi_i\cap\pi_j\cap\pi_k$ with the line $L_{lm} = \pi_l\cap
\pi_m$, where $\{i,j,k,l,m\}=\{1,2,3,4,5\}$.
 The involution $\tau$ of $\tilde H$ is fixed-point-free
and its quotient $Y = \tilde H/(\tau)$ is an Enriques surface. The
pair of points in the involution defines a line in $\bbP^3$. The set of
such lines forms the Reye congruence of lines isomorphic
to
$Y$ (see [Co]).
\item{(8)} The pencil of planes passing through a line $L_{ij}$ of
$H$ cuts out a pencil of cubic curves on $H$. The plane touching
$H$ along the line $L_{ij}$ defines a conic $C_{ij}$ on $H$.
\item{(9)} The union of the conics corresponding to the four lines lying
in
the same plane $\pi_i$ is cut out by a quadric.
\item{(10)} Projecting from a node $P_{ijk}$  we get a representation
of
$H$ as the double cover of $\bbP^2$ branched along the union
of two cubics. The two cubics intersect at the six vertices of a complete
quatralaterial $xyz(x+y+z) = 0$, tangent at three vertices
$(0,0,1),(0,1,0),(1,0,0)$ with the tangent directions such that there
exists a conic intersecting each cubic only at these three points.

\head{\bf 2. The Reye congruence}\endhead

\bigskip\noindent
The Hessian surface $H$ is a special case of a quartic symmetroid, that
is, a quartic surface given by vanishing of the determinant of a
symmetric matrix with linear homogeneous polynomials as its entries.  We
refer to [Co] for  general properties of quartic symmetroids and the
associated Reye congruences of lines.

Let us look closer at the Enriques surface $Y = \tilde H/(\tau)$.
We have a natural indexing of faces $\pi_i$ of the Sylvester
pentahedron by the set $\{1,2,3,4,5\}$, edges by subsets $\alpha
= \{i,j\}\subset \{1,2,3,4,5\}$ and the vertices by subsets
$\beta = \{i,j,k\}\subset \{1,2,3,4,5\}$. A vertex $P_{ijk}$
belongs to an edge $L_{mn}$ if and only if $\{i,j,k\}\supset
\{m,n\}$. Let $N_\alpha$ be the smooth rational curve on $\tilde
H$ corresponding to $P_\alpha$ and similarly let $T_\beta$ be the
smooth rational curve corresponding to $L_\beta$. The Enriques
involution $\tau$ interchanges the curves $N_\alpha$ and
$T_\beta$, where $\alpha$ and $\beta$ are complementary sets. Let
$U_{\beta}$ be the image on $Y$ of the pair $(N_\alpha,T_\beta)$.
We have
$$U_{ab}\cdot U_{cd} = \cases 1,&\text{if $\{c,d\}\cap \{a,b\} =
\emptyset$};\\
0,&\text{otherwise.}\endcases.$$
The dual intersection graph is the famous Peterson tri-valent graph
with group of symmetry isomorphic to the permutation group
$S_5$.

\epsfxsize = 200pt \centerline{\epsfbox{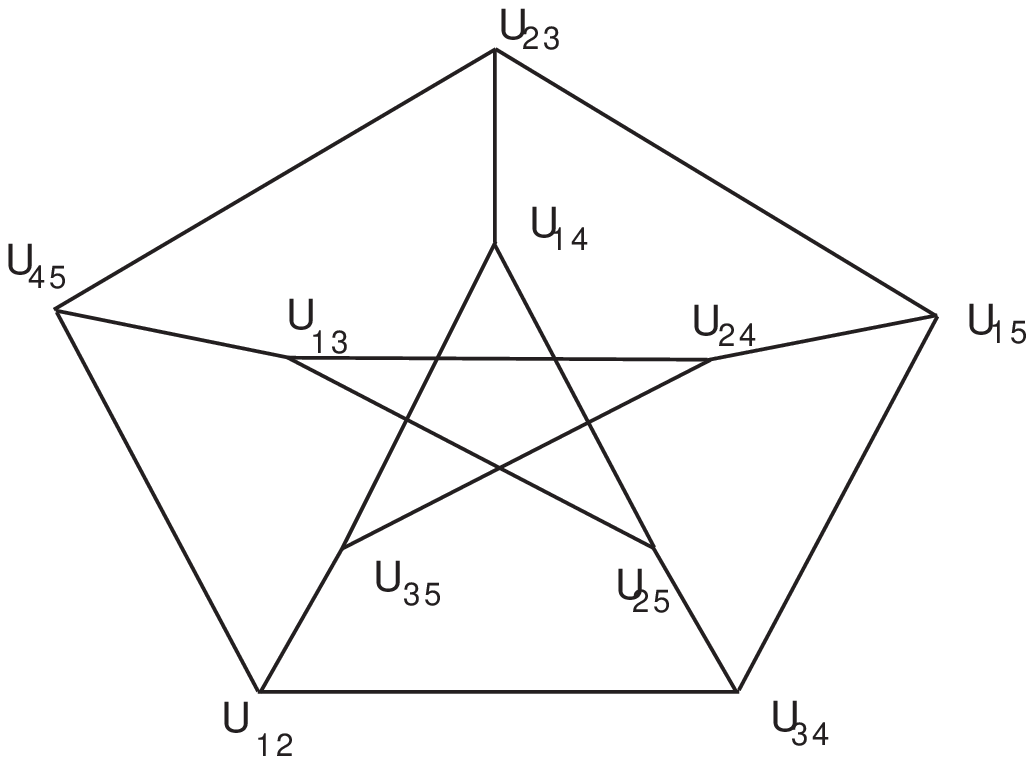}}
\bigskip
\centerline{Fig.1}
\bigskip
The sum $\Delta$ of the nodal curves $U_{ab}$ satisfies $\Delta^2
= 10$. It defines an embedding of $Y$ into the Grassman variety
$G(2,4)$ which exhibits $Y$ as the Reye congruence of lines for
the net of polar quadrics to $F$.

Let $D_{12} = U_{12}+U_{34}+U_{35}+U_{45}$ and $F_{12} =
\Delta-D_{12}$. Then $|F_{12}|$ is an elliptic pencil on $Y$ with
degenerate fibre $F_{12}$ of Kodaira type $I_6$ and one fibre of type
$I_2$ formed by $U_{12}$ and the image of the conic $C_{12}$.
Similarly we define other elliptic pencils $|F_{ab}|$. There are
ten such elliptic pencils on $Y$. We have $F_{ab}\cdot F_{cd} = 4$
if $\{c,d\}\ne \{a,b\}$. Note that
$$\sum_{ab}F_{ab}=6\Delta,\quad \sum_{ab}E_{ab}=3\Delta,$$
where $2E_{ab}$ is a multiple fibre of $|F_{ab}|$. If $U_{ab}$ and
$U_{cd}$ intersect, the sum $D_{ab}+D_{cd}$ is a degenerate fibre
of type $I_2^* = \tilde D_5$ of the elliptic pencil
$|2\Delta-F_{ab}-F_{cd}|=|D_{ab}+D_{cd}|$. The linear system
$|F_{ab}+F_{cd}|$ defines a finite map of degree 2 from $Y$ to a
Del Pezzo surface $\calD_4$ of degree 4. It blows down two Dynkin
curves of type $A_2$. The ramification curve is of genus 3.

The pull-backs of the elliptic pencils on $Y$ to $\tilde H$ give
elliptic pencils of two types: Type I with degenerate fibres
$2I_6+2I_2$ and Type II with degenerate fibres $2I_2^*$. A
pencil of Type I is cut out by planes through a line $L_{ij}$ or, in the
double plane construction, by the pencil of cubic curves spanned by
the components of the branch locus.

\head{\bf 3. The Picard lattice}\endhead

\bigskip\noindent
Let us consider the double plane construction of $H$ and the
corresponding elliptic pencil. The surface $\tilde H$ is the
double cover of a rational elliptic surface $V$ obtained by
blowing up the base points of the pencil of cubic curves spanned
by the components of the branch locus. It is branched along two
elliptic fibres. The Mordell-Weil group of $V$ is of rank 2 and
is isomorphic to the Mordell-Weil group of the elliptic pencil on
$\tilde H$ (obtained from the pencil on $V$ by the double cover
of the base curve ramified at two points). This implies that any
divisor on $\tilde H$ can be written uniquely as a sum $D_1+D_2$,
where $D_1$ is the pre-image of a divisor on $V$ and $D_2$ is a
combination of irreducible components $E_i, i = 0,\ldots,5,$ of a
fibre of type  $I_6$ and one component $E_6$ of a fibre of type
$I_2$. Consider the sublattice $S$ of $\Pic(\tilde H)$ spanned by
$\pi^*(\Pic(V))$ and the divisors $D_i = E_i-\alpha(E_i)$, where
$\pi:\tilde H\to V$ is the covering map and $\alpha:\tilde H\to
\tilde H$ the covering involution. We have $\sum_{i=0}^5 D_i = 0$
and $2E_i\in S$. This gives
$$S \cong \pi^*(\Pic(V))\perp \sum_{i=1}^6 \bbZ[D_i] \cong U(2)\perp
E_8(-2)\perp A_5(-2)\perp A_1(-2),$$
$$\Pic(\tilde H)/S \cong (\bbZ/2)^7.$$
Here $U$ denotes the standard hyperbolic plane, and $M(m)$ denotes the
lattice $M$ with quadratic form multiplied by $m$. The generators of the
latter group are the classes of the divisors
$E_i, i = 0,\ldots,6$. This shows that the discriminant of the lattice
$\Pic(\tilde H)$ is equal to
$2^2\cdot 2^8\cdot 2^5\cdot 6\cdot 2\cdot 2/2^{14} = 2^4\cdot 3.$ More
lattice theoretical computation shows that the discriminant quadratic
form of $\Pic(\tilde H)$ coincides with the discriminant quadratic
form of the lattice $U\perp U(2)\perp A_2(-2)$. Applying
Nikulin's results [Ni] we obtain that the transcendental lattice
of $\tilde H$ is isomorphic to
$$T_{\tilde H} \cong U\perp
U(2)\perp A_2(-2).$$

\medskip
\plainproclaim Lemma 3.1. The Picard lattice $\Pic(\tilde H)$ is
generated over integers by the twenty smooth rational curves
$N_{\alpha}$ and $T_{\beta}$.

{\sl Proof.} We have seen already that $\Pic(\tilde H)$ is generated by
$\pi^*(\Pic(V))$ and the curves $E_i, \alpha(E_i), i = 0,\ldots,5$. The
latter curves belong to the set $NT = \{N_\alpha,T_\beta\}$. The
group $\Pic(V)$ is generated by the exceptional curves blown up
from the vertices of the quatralaterial and the proper
transforms of its sides. The pre-images of all these curves in
$\tilde H$ belong to the set $NT$.

\medskip
 We set
$$\calN = \sum_\alpha N_\alpha, \quad \calT = \sum_\beta T_\beta.$$
Let $\eta_H$ be the pre-image on $\tilde H$ of the class of a
hyperplane section of $H$ and let $\eta_S$ be its image under the
Enriques involuton $\tau$. It is known (see, for example, [Co],
Proposition 2.4.1) that
$$2\eta_H = 3\eta_S-\calT,$$
$$2\eta_S = 3\eta_H-\calN.\eqno (3.1)$$
In particular, the pre-image $\tilde \Delta$ of the class
$\Delta$ on the Enriques surface can be expressed as
$$\tilde\Delta = \calT+\calN = \eta_H+\eta_S.\eqno (3.2)$$
Pick up a face of the Sylvester pentahedron, say $\pi_5$. It is
immediately checked that $2(\eta_S-\eta_H)$ intersects each $N_\alpha$ and
$T_\beta$ with the same multiplicity as the divisor
$(T_{15}+T_{25}+T_{35}+T_{45})-(N_{123}+N_{124}+N_{134}+N_{234}).$
This gives a linear relation:
$$2(\eta_S-\eta_H) =
(T_{15}+T_{25}+T_{35}+T_{45})-(N_{123}+N_{124}+N_{134}+N_{234}).$$
Denote the first bracket by $\calT_{5}$ and the second bracket by
$\calN_{5}$. Similarly introduce $\calT_i, \calN_i, i =
1,\ldots,4,$ for any other face. We obtain
$$2(\eta_S-\eta_H) = \calT_i-\calN_i, \quad i = 1,\ldots, 5.\eqno (3.3)$$
One can show that the linear system $|\eta_H+{1\over
2}(\calT_i-\calN_i)| = |\eta_S|$ defines a standard cubic Cremona
involution $\tau_i$ with fundamental points at the vertices of
the tetrahedron formed by the faces $\pi_j, j\ne i$. All these
involutions restrict to the Enriques involution on $\tilde H$ (cf.
[Hu2], p.335).

We shall denote by $C_{ij}$ the pre-image on $\tilde H$ of the
conic cut out by the plane which is tangent to $H$ along the edge
$L_{ij}$. Its divisor class equals
$$C_{ij} = \eta_H-2T_{ij}-N_{ijk}-N_{ijl}-N_{ijm}.\eqno (3.4)$$
We shall also denote by $R_{klm}$ the pre-image on $\tilde H$ of
the cubic cut out by the plane through the edge $L_{ij}$ and the
opposite vertex $P_{klm}$. Its divisor class is equal to
$$R_{klm} = \eta_H-T_{ij}-N_{ijk}-N_{ijl}-N_{ijm}-N_{klm}.\eqno (3.5)$$
One can show ([Hu2]) that any rational nonsingular curve of degree $\le
3$ is one of the described above (on a general $H$).

\head {\bf 4. Kummer surfaces and Hessians}\endhead

\bigskip\noindent

Let $C$ be a genus 2 curve and  $(p_1,\ldots,p_6)$ be its ordered
set of Weierstrass points. Then the divisor classes $\mu_0 = 0,\mu_{ij} =
[p_i+p_j-2p_6], 1\le i< j\le 6,$ are the sixteen 2-torsion points on its
Jacobian variety $J$. In the usual way we identify the group of 2-torsion
points $J_2$ with the set of 2-element subsets
$\alpha$ of the set
$S =
\{1,2,3,4,5,6\}$, the zero point  corresponds to the
empty subset.
 We shall identify a 2-element subset $\alpha$ with its
complementary subset $S\setminus \alpha$. Then the addition of points
corresponds to the symmetric sum $\alpha+\beta$ of subsets. One defines
the symplectic bilinear form with values in $\bbF_2$ on
$J_2$ by
$$(\mu_\alpha,\mu_\beta) = |\alpha\cap\beta|\quad \roman{modulo}\quad
2.\eqno (4.1)$$ Let us write elements of the vector space $V =
\bbF_2^4$ as $(2\times 2)$-matrices $[\epsilon,\eta]$ with
columns in $\bbF_2^2$. We define the quadratic form $q_0:V\to
\bbF_2$  by $q_0([\epsilon,\eta]) = \epsilon\cdot \eta,$ where
the dot means the standard dot-product in $\bbF_2^2$. The
associated symmetric bilinear form is
$$([\epsilon,\eta],[\epsilon',\eta']) = \epsilon\cdot \eta'+\eta\cdot \epsilon'.$$
It is a non-degenerate symplectic form. We define an isomorphism
of symplectic spaces $\psi:J_2\to V$ by
$$\mu_{12}\to \left[\matrix 1&0\\
0&0\endmatrix\right],\quad \mu_{34}\to \left[\matrix 0&1\\
0&0\endmatrix\right],\quad \mu_{16}\to \left[\matrix 0&0\\
1&0\endmatrix\right],\quad \mu_{45}\to \left[\matrix 0&0\\
0&1\endmatrix\right].$$ Let us identify the span $U_1$ of
$\mu_{12},\mu_{16}$ with the set $[4] =\{1,2,3,4\}$ by assigning
$0$ to $1$, $\mu_{16}$ to 2, $\mu_{12}$ to 3 and
$\mu_{12}+\mu_{16} = \mu_{26}$ to $4$. Similarly we identify the
span $U_2$ of $\mu_{34},\mu_{45}$ with $[4]$ by assigning $0$
to $1$, $\mu_{45}$ to 2, $\mu_{34}$ to 3 and $\mu_{34}+\mu_{45} =
\mu_{35}$ to $4$.  Each 2-torsion point can be written uniquely
as the sum $a+b, a\in U_1,b\in U_2$, and hence can be identified
with the pair $(a,b)\in [4]\times [4]$, or, equivalently, with
a dot in the following $(4\times 4)$-table:
$$\matrix\bullet&\bullet&\bullet&\bullet\\
\bullet&\bullet&\bullet&\bullet\\
\bullet&\bullet&\bullet&\bullet\\
\bullet&\bullet&\bullet&\bullet\endmatrix$$

 Let $C\to J, x\to [x-p_6],$ be the Abel-Jacobi map corresponding
to the point $p_6$. The image of $C$ is denoted by $\Theta_0$. Let
$\Theta_{ij} = \Theta_0+\mu_{ij}$ be its translate by a 2-torsion
point $\mu_{ij}$. Each $\Theta_{ij}$ contains exactly six
2-torsion points, namely the points $\mu_{ij},\mu_{k6}+\mu_{ij},
k = 1,\ldots, 5$. In other words,
$$\mu_{\alpha}\in \Theta_\beta \Leftrightarrow \beta+\alpha \in
\{\emptyset, \{16\}, \{26\}, \{36\}, \{46\}, \{56\} \}.$$ One also
employs different indexing of theta divisors $\Theta_{\alpha}$.
For each $\alpha$ there exists a partition of $[6] = \{1,2,3,4,5,6\}$ into
two disjoint subsets $S\cup S'$ with odd number elements. For
$\alpha\neq k6$, it is uniquely determined by the property
$\mu_\beta\in \Theta_\alpha$ if and only if $\beta\subset S$ or
$\beta\subset S'$. For $\alpha=k6$, we simply take $S=\{k\}$. We
use either $S$ or $S'$ for the index. In this correspondence,
$\Theta_{\alpha}$ corresponds to $\Theta_{\alpha+\{6\}}$. For
example, $\Theta_{12}$ corresponds to $\Theta_{345}$ or
$\Theta_{126}$, and $\Theta_{16}$ to $\Theta_{1}$ or
$\Theta_{23456}$. Yet there is another classical notation for
a theta divisor. Each theta divisor is equal to the set of zeroes
of a theta function
$\theta\left[\matrix\epsilon\\
\epsilon'\endmatrix\right](z,\tau)$ with theta characteristic $\left[\matrix\epsilon\\
\epsilon'\endmatrix\right]$. The theta characteristic
corresponding to $\Theta_{S},$ where $\#S $ is odd, is equal to
$\psi(S+\{135\})$. For example, the theta characteristic
corresponding to $\Theta_{12}$ is
$\psi(14)=\left[\matrix1&1\\
1&1\endmatrix\right]$, where $\epsilon$ and $\epsilon'$ are the
first and the second columns of the matrix.

 One assigns to
$\Theta_\alpha$ a dot $(ab)\in [4]\times [4]$ in the
$(4\times 4)$-table as above in such a way that $(ab)$ in the
right-side table contains a 2-torsion point corresponding to the
entry $(cd)$ in the left-side table if and only if $a = c$ or $b
= d$ but $(ab)\ne (cd)$.

Explicitly, we have two tables
$$\matrix\mu_0&\mu_{45}&\mu_{34}&\mu_{35}\\
\mu_{16}&\mu_{23}&\mu_{25}&\mu_{24}\\
 \mu_{12}&\mu_{36}&\mu_{56}&\mu_{46}\\
\mu_{26}&\mu_{13}&\mu_{15}&\mu_{14}\endmatrix \qquad\qquad
\matrix\Theta_{12}&\Theta_{36}&\Theta_{56}&\Theta_{46}\\
\Theta_{26}&\Theta_{13}&\Theta_{15}&\Theta_{14}\\
\Theta_{0}&\Theta_{45}&\Theta_{34}&\Theta_{35}\\
\Theta_{16}&\Theta_{23}&\Theta_{25}&\Theta_{24}\endmatrix$$

The following construction of Hutchinson [Hu1] describes the
translation of left-hand-side table to the table of the
corresponding values of the map $\psi$, and, at the same time,
translates the right-hand-side table to the table of the
corresponding theta characteristics. Denote the four columns
$$\matrix 1&1&0&0\\
1&0&1&0\endmatrix$$ by the numbers $1,2,3,4$, respectively. Then
the theta characteristic of the $(ab)$-entry in the
right-hand-side table is equal to the characteristic formed by the columns
$(ab)$. Now, if reverse the order of the 4 vectors and do the
same, we obtain the values of $\psi$ at the entries of the
left-hand-side table.

\bigskip
 An affine plane in
$J_2$ is called an odd (or G\"opel) tetrad if it is a translation
of a totally isotropic linear subspace of dimension 2. Otherwise
it is called an even (or Rosenhain) tetrad.  There are 60 odd and
80 even tetrads, each set forming an orbit with respect to the
group generated by symplectic automorphisms of $J_2$ and
translations. For example, the rows and columns of the left-hand-side
table above correspond to even tetrads but the diagonal corresponds to an
odd tetrad.

A {\it Weber hexad} is defined as the symmetric sum of an even
and an odd tetrad which have one point in common.  For example,
the dots marked with asterisk represent a Weber hexad:
$$\matrix\ast&\bullet&\ast&\bullet\\
\bullet&\ast&\ast&\bullet\\
 \bullet&\bullet&\bullet&\bullet\\
\bullet&\bullet&\ast&\ast\endmatrix.$$ A Weber hexad has the
following property. Each theta-divisor contains either 3 points
from a Weber hexad or just one. The number of theta divisors
which contain three points is equal to ten. One can choose 5
subsets of four elements from this set of ten theta divisors such
that each point from the Weber hexad is contained in two theta
divisors from this set.  For the Weber hexad chosen in above, the
five sets are the following:
$$\calA_1 = (\Theta_{56},\Theta_{46},\Theta_{15},\Theta_{14}),\quad
\calA_2 = (\Theta_{14},\Theta_{36},\Theta_{16},\Theta_{34}),$$
$$\calA_3 = (\Theta_{23},\Theta_{25},\Theta_{56},\Theta_{36}),\quad
\calA_4 = (\Theta_{23},\Theta_{26},\Theta_{34},\Theta_{46}),$$
$$\calA_5 = (\Theta_{26},\Theta_{16},\Theta_{25},\Theta_{15}).\eqno (4.2)$$
By inspection, one sees that theta  characteristics of theta
functions from the same group of four add up to zero. This means
that the sum of the divisors belong to the linear system
$|4\Theta|$, where $\Theta = \Theta_{24}$ has the zero theta
characteristic. Thus we have five divisors $D_1,\ldots,D_5\in
|4\Theta|$ each passing through every point of the Weber hexad
with multiplicity 2. Let $\theta_i, i = 0,\ldots,4,$ be the
corresponding theta functions of order 4 with zero theta
characteristic. They define a rational map $\Phi:J\to \bbP^4$.

\plainproclaim Theorem 4.1. The image $H$ of $\Phi$ is contained
in a linear hyperplane $x_4 = ax_0+bx_1+cx_2+dx_3,$ where
$a,b,c,d\ne 0$, and satisfies the quartic equation:
$$ax_1x_2x_3x_4+bx_0x_2x_3x_4+cx_0x_1x_3x_4+dx_0x_1x_2x_4+x_0x_1x_2x_3
= 0.$$ In particular, $H$ is isomorphic to the Hessian of a nonsingular
cubic surface.

\medskip
We refer for the proof to [Hu1] or [vGe]. One can also find there
the expression of the coefficients $a,b,c,d$ in terms of theta
constants.

\medskip
Since $\Theta$ is a symmetric theta divisor, the rational map
$\Phi$ factors through the map $J \to \Kum(J)\subset \bbP^3$ given
by the linear system $|2\Theta|$. The image of the Weber hexad in
the Kummer surface is a set of 6 nodes, a Weber set of nodes. The
linear system $|4\Theta|$ is equal to the inverse transform of a
linear system $|L_h|$ of quadrics through the Weber set of nodes
$h$.  Thus the Hessian surface $H = \Ima(\Phi)$ is the birational
image of $\Kum(J)$ under the map $\bar \Phi$ given by the linear
system $|L_h|$.   An explicit equation of the locus of the cubic
surfaces whose Hessian is birationally isomorphic to a Kummer
surface inside of the moduli space of cubic surfaces has been
found in [Ro,vGe].

Let us choose the Weber hexad as above:
$$h = (\mu_0,\mu_{23},\mu_{34},\mu_{25},\mu_{15},\mu_{14}).\eqno (4.3)$$
Let
$$(\Theta_{56},\Theta_{46},\Theta_{15},\Theta_{14},
\Theta_{36},\Theta_{16},\Theta_{34},\Theta_{23},\Theta_{25},\Theta_{26})
\eqno (4.4)
$$ be the corresponding set of ten theta divisors which contain exactly
three of the points from $h$ and hence three points from the set
$J_2\setminus h$.  The images of each node from $h$ under the map
$\Phi$ is a conic in $H$. The plane section of $H$ along such a
conic
 is equal to the union of two conics. These conics do not appear on a
general Hessian.

The images of the remaining nodes which are the images of ten
2-torsion points $\mu_\alpha$ are ten nodes of the Hessian which
we denote by $N_\alpha$. The image of a theta divisor
$\Theta_\beta$ from (4.4) is a line on $H$ which contains three
nodes $N_\alpha$ such that $\mu_\alpha\in \Theta_\beta$. We
denote this line by $T_{\beta}$. Note that each $\Theta_\beta$
belongs to exactly two of the five subsets from (4.2), say
$\calA_i, \calA_j$.  If  we reindex $T_\beta$ by $T_{ij}$ and
$N_\alpha = N_{ijk}$, where $N_\alpha\in T_{ij}$, $T_{jk}$,
$T_{ik}$, then we get the notation used for the nodes and the
lines on a general Hessian. The image of a theta divisor not
belonging to the group of ten is a rational curve of degree 3 on
$H$ which passes through 5 nodes. These curves do not appear on a
general Hessian.

 \medskip\noindent
{\bf Remark 4.2}. A Weber set $h$ of nodes of a Jacobian Kummer
surface $K$ is a set $\calP$ of 6 points in $\bbP^3$ in general
linear position.  Recall that each such set defines the Weddle
quartic surface $W(\calP)$ which is  the locus of singular points of
quadrics passing through the set $\calP$. The image $W(\calP)$ by
the linear system of quadrics through $\calP$ is a Kummer surface
$K(\calP)$ birationally isomorphic to $W(\calP)$. The image of $K$ is a
Hessian quartic
$H$. The  quartic surfaces $K(\calP)$ and $H$
 touch each other  along a curve of degree $8$ not passing
through their nodes. There are 192 Weber hexads. The affine symplectic
group $2^4\rtimes \roman{Sp}(4,\bbF_2)$ of order
$2^4\cdot 6$! acts transitively on the set of Weber hexads. The
isotropy group of a Weber hexad is isomorphic to the alternating
group $A_5$. Let $\calM_{\Kum}^w$ be the moduli space of Jacobian
Kummer surfaces together with a choice of a Weber hexad of nodes.
It is a
cover of degree 12 of the moduli space $\calM_2$ of genus 2
curves isomorphic to the moduli space $\calM_{\Kum}$ of Jacobian
Kummer surfaces. The construction from above defines a map
$\calM_{\Kum}^w\to \calM_{\Kum}$. What is the degree of this map?

 \medskip\noindent
{\bf Remark 4.3}. Since the K3-surface birationally isomorphic to a
general Hessian quartic admits an Enriques involution, any
K3-surface birationally isomorphic to a Jacobian Kummer surface
admits an Enriques involution. This fact, of course is known, and
also true for not necessarily Jacobian Kummer surfaces (see
[Ke1]). But here we get an explicit construction of this
involution.

Applying (3.1), we see that the Enriques involution on the
Hessian quartic surface associated to a Kummer surface $K$ is
given by the linear system defined by the divisor
$$\eta_S \sim {1\over 2}(3\eta_H-\calN) \sim {1\over
2}(3(2\eta_K-\sum_{p\in W}E_p)-\sum_{p\not\in W}E_p) =
3\eta_K-{3\over 2}\sum_{p\in W}E_p-{1\over 2}\sum_{p\not\in
W}E_p,$$ where $\eta_K$ is the class of a hyperplane section of
the Kummer surface $K$ in $\bbP^3$, $W$ is the Weber hexad of
nodes corresponding to the Weber hexad of 2-torsion points, and
$E_p$ is the class of the exceptional curve blown up from one of
the 16 nodes of $K$. On the abelian surface $J$ this corresponds
to the linear system
$$|6\Theta-3\sum_{\mu\in h}\mu-\sum_{\mu\not\in
h}\mu|.$$ The linear system
$$|6\Theta-3\sum_{\mu\in h}\mu-\sum_{\mu\not\in
h}\mu+(4\Theta-2\sum_{\mu\in h}\mu)| =|
10\Theta-5\sum_{\mu\in h}\mu- \sum_{\mu\not\in
h}\mu|$$ maps $J$  onto an Enriques surface embedded in
$\bbP^5$ by its Fano linear system of degree 10.

If we choose the Weber hexad (4.3),
then the Sylvester pentahedron of the corresponding Hessian quartic is given
by the following figure:

\bigskip
\epsfxsize = 250pt

\centerline{\epsfbox{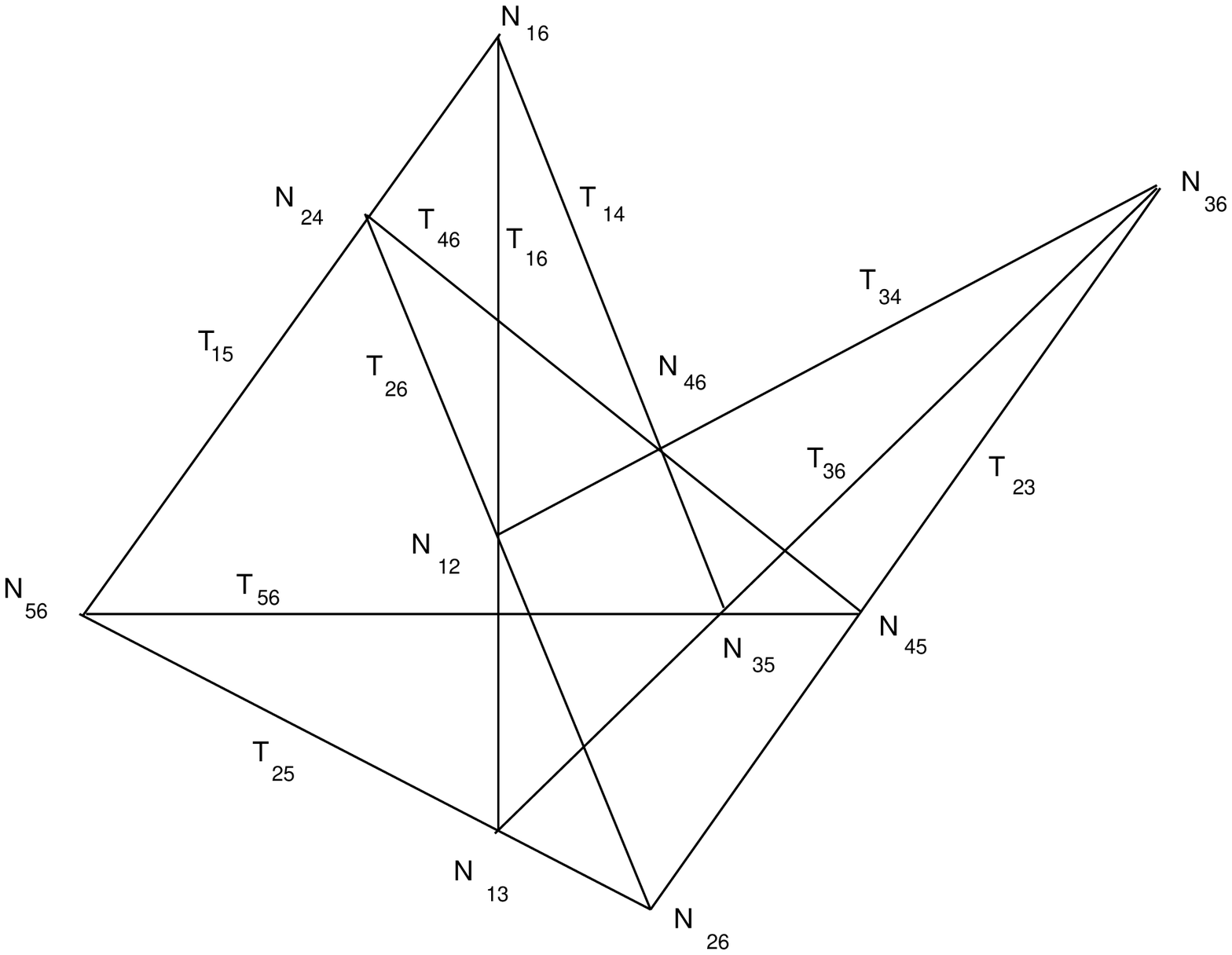}}
\centerline{Fig.2}

We shall  use the previous notation for nodes and lines on a
general Hessian surface.

\head {\bf 5. Elliptic pencils}\endhead
\bigskip\noindent
There are many elliptic pencils on the surface $\tilde H$. We
list a few, which will be used later. We accomodate the notations
for the conics $C_\alpha$ (3.4) and cubics $R_\alpha$ (3.5) to
our new notation for (-2)-curves $T_\alpha$ and $N_\beta$.

\bigskip\noindent
{\bf Type 1.} It is cut out by the linear pencil of planes
through an edge of the pentahedron:
$$|F_{\alpha}| = |C_\alpha+T_\alpha|.$$
We have
$$ F_{\alpha}\cdot F_{\beta} = 2, \quad \hbox{ if $T_{\alpha}$ and $T_{\beta}$ are skew}.$$
The
pencil
$|F_\alpha|$ has 2 reducible fibres of type
$I_2$:
$$C_\alpha+T_\alpha, \quad N_{\bar \alpha}+R_{\bar\alpha},$$
where $R_{\bar\alpha}$ is the residual cubic for the plane section
of $H$ passing through the edge $T_\alpha$ and its opposite node
$N_{\bar\alpha}$ (see (3.5)). It has also two reducible fibres of
type $I_6$ corresponding to the two faces containing $T_\alpha$.

Observe that the Enriques involution leaves the pencil invariant by
interchanging two degenerate fibres of the same type. The members of the
pencil are cubic curves.

The pencil of type 1 is denoted by $F_i$ in [CD]. It is the
pre-image of one of the ten pencils on the Enriques surface $Y$
as explained in Section 2. Using $(3.2)$, we have
$$|\tilde \Delta-2F_i| = B_i+\tau(B_i)$$
where $\tau$ is the Enriques involution, and $B_i$ is of the form
$T_{15}+N_{56}+N_{24}+N_{16}$. This agrees with the notation from
[CD].

\bigskip\noindent
{\bf Type 2.} This is a pencil $|F_{N_\alpha,T_\beta}|$ on $\tilde H$
formed by proper transform of quartic elliptic curves  cut out by the
pencil of quadric cones with the vertex at a node $N_\alpha$ which
contain the lines  through
$N_\alpha$ and tangent along one of them, $T_\beta$.   For example
$$|F_{N_{16},T_{15}}| =
|2\eta_H-2T_{15}-T_{16}-T_{14}-N_{24}-N_{56}-N_{12}-N_{13}-N_{46}-N_{35}-2N_{16}|$$
It has one reducible fibre of type $I_8$:
$$T_{25}+N_{26}+T_{26}+N_{24}+T_{46}+N_{45}+T_{56}+N_{56},$$
and two reducible fibres of type $I_4$:
$$C_{15}+N_{36}+T_{34}+T_{36},\quad C_{23}+N_{16}+T_{16}+T_{14}.$$
Observe that $N_{46},N_{35},N_{13},N_{12}$ are sections and $T_{23},
T_{15}$ are bisections of the elliptic fibration. Let us take $N_{46}$ as
the zero section. We easily check (by intersecting both sides with any
$N_\alpha$ and any $T_\beta$) that
$$2N_{13}-
2N_{46}\sim
-T_{36}+T_{34}+T_{14}-T_{16}+T_{46}-N_{56}-2T_{25}-N_{26}+N_{24}+2T_{46}+N_{45}.$$
This implies that $N_{13}$ is a 2-torsion section. Also
$$N_{46}+N_{12} \sim
N_{13}+N_{35}-N_{46}-T_{34}+T_{36}-T_{46}+T_{56}+N_{56}+T_{25}-T_{26}-N_{24}$$
implies that $N_{13}\oplus N_{35} = N_{12}.$ One also checks that
translation by $N_{13}$ sends $T_{15}$ to $R_{36}$ and $T_{23}$ to
$R_{16}$.

Note that
$$|F_{N_{16},T_{15}}|=|F_{N_{36},T_{23}}|.$$
There are 15 elliptic pencils of type 2.


\bigskip\noindent
{\bf Type 3.} This is a pencil $|F_{T_\beta,T_\gamma}|$ on $\tilde H$ of
proper transforms of quartic elliptic curves  which are cut out by the
pencil of quadric cones with the vertex at a node $N_\alpha$ and tangent
to $H$ along two edges $T_\beta,T_\gamma$ intersecting at $N_\alpha$. It
is spanned by the double plane through the two edges and the union of two
planes tangent along these edges. For example
$$|F_{T_{16},T_{14}}| =
|2\eta_H-2T_{16}-2T_{14}-N_{12}-N_{13}-N_{46}-N_{35}-2N_{16}|=
|C_{16}+C_{14}|.$$ It has  reducible fibre of type $ I_0^*$ and
$I_2^*$:
$$C_{23}+2T_{15}+N_{16}+N_{56}+N_{24},$$
$$2T_{36}+2T_{34}+2N_{36}+N_{13}+N_{12}+N_{46}+N_{35}.$$
Also it has a reducible fibre of type $I_2$:
$$C_{16}+C_{14}.$$
We see that $N_{45}$ and $N_{26}$ are contained in fibres. Since
$C_{16}+C_{14}-N_{26}$ (resp. $C_{16}+C_{14}-N_{45}$) is of degree
4 with respect to $\eta_H$, and cannot be a union of two conics,
we get two more reducible fibres of type $I_2$. Observe that
$T_{26},T_{56},T_{25}$ and $T_{46}$ are sections and
$T_{16},T_{14},T_{23}$ are bisections. Let us take $T_{46}$ as
the zero section. Then
$$2T_{25}-2T_{46} \sim
2T_{36}+3N_{36}+4T_{34}+N_{35}+3N_{46}+2N_{12}-N_{56}+$$
$$N_{24}+C_{16}+N_{45}-N_{26}-2(C_{16}+C_{14})$$
implies that $T_{25}$ is a 2-torsion section. Also
 $T_{25}\oplus T_{56} = T_{26}.$


\head {\bf 6. Birational involutions of a Hessian quartic}\endhead

There are some obvious birational involutions of a Hessian
quartic surface. They are:

\bigskip
\item{(i)} the Enriques involution $\tau$;
\item{(ii)} an involution $p_{\alpha}$ defined by projection from a node $N_{\alpha}$;
\item{(iii)} an involution defined by a pair of skew lines;
\item{(iv)} an involution defined by two elliptic pencils $|F|$ and $|F'|$
with $F\cdot F' = 2$.  This is the covering involution for a
degree 2 map $\tilde H \to \bbP^1\times \bbP^1$ defined by the
linear system, $|F+F'|$.
\item{(v)} an involution defined by the inversion map of an elliptic
fibration with a section;
\item{(vi)} an involution defined by
the translation by a 2-torsion section  in the group law of
sections of an elliptic fibration.

\bigskip
Let us describe the action of each involution on the Picard
lattice $S_H$ of $\tilde H$.

\bigskip
(i) We have already described the action of the Enriques involution.
In our new notation we have:
$$N_{16}\longleftrightarrow T_{23}, N_{24}\longleftrightarrow T_{36},
N_{56}\longleftrightarrow T_{34}, N_{12} \longleftrightarrow T_{56},
N_{13} \longleftrightarrow T_{46}$$
$$N_{26} \longleftrightarrow T_{14}, N_{35} \longleftrightarrow T_{26},
N_{46}\longleftrightarrow T_{25}, N_{36}\longleftrightarrow T_{15},
N_{45}\longleftrightarrow T_{16}.$$
 We already know that the image of $\eta_H$ is equal to $\eta_S$.

\bigskip
(ii) We may assume that we project to the plane defined by the
opposite face of the pentahedron. Then the projection $p_{16}$
from $N_{16}$ acts as follows:
$$N_{56}\longleftrightarrow N_{24}, N_{13}\longleftrightarrow N_{12},
N_{35}\longleftrightarrow N_{46},$$
$$T_{25}\longleftrightarrow T_{26},
T_{36}\longleftrightarrow T_{34}, T_{56}\longleftrightarrow T_{46},$$
$$T_{23}\longleftrightarrow \eta_H-N_{16}-T_{23}-N_{26}-N_{45}-N_{36} =
R_{16},$$
and the remaining curves $N_{26}, N_{36}, N_{45}, T_{15}, T_{16},
T_{14}$ are fixed.

By (3.4-3.5),  $N_{16}+R_{16} = C_{23}+T_{23}$, hence
$$N_{16} \longleftrightarrow C_{23}.$$
Since $\eta_H = C_{23}+2T_{23}+N_{36}+N_{45}+N_{26},$ we obtain
$$\eta_H \longleftrightarrow N_{16}+2R_{16}+(N_{36}+N_{45}+N_{26}) =
2\eta_H-2T_{23}-N_{26}-N_{45}-N_{36}-N_{16}.$$ Thus the involution is
given by the linear system of quadrics through the vertex $N_{16}$
and touching $H$ along the edge $T_{23}$.

Projections $p_{\alpha}$ commute with the Enriques involution
$\tau$, i.e. $$p_{\alpha}\circ\tau=\tau\circ p_{\alpha}.$$

\bigskip (iii) This is a special case of (iv) when the two pencils
are of type 1 and correspond to two skew edges. Assume that the
edges are $T_{15}$ and $T_{23}$. Let the fibres
$N_{12}+T_{16}+N_{13}+T_{25}+N_{26}+T_{26}$ and
$T_{14}+T_{46}+N_{46}+N_{15}+T_{56}+N_{35}$ of $F_{15}$ go to the
lines $A_1,A_2$ on the quadric $\bbP^1\times \bbP^1$. Similarly,
let $T_{56}+T_{36}+T_{25}+N_{56}+N_{35}+N_{13}$ and
$T_{26}+T_{46}+T_{34}+N_{12}+N_{46}+N_{24}$ go to the lines
$B_1,B_2$ of the other ruling. Then the pre-image of $A_1$ splits
in $T_{16}+N_{26}$, the pre-image of $A_2$ splits in
$T_{14}+N_{45}$, the pre-image of $B_1$ splits in
$T_{36}+N_{56}$, the pre-image of $B_2$ splits in
$T_{34}+N_{24}$. This easily shows that the action of the
involution is defined as follows:
$$T_{16}\longleftrightarrow N_{26}, T_{14}\longleftrightarrow N_{45},
T_{36}\longleftrightarrow N_{56}, T_{34}\longleftrightarrow N_{24},
$$
$$T_{15}\longleftrightarrow C_{15}, N_{36}\longleftrightarrow R_{36},
N_{16}\longleftrightarrow R_{16}, T_{23}\longleftrightarrow
C_{23},$$
$$T_{25}\longleftrightarrow N_{13}, T_{56}\longleftrightarrow N_{35},
T_{26}\longleftrightarrow N_{12}, T_{46}\longleftrightarrow N_{46}.$$ This
implies that
$$\eta_H\longleftrightarrow  T_{15}+2C_{15}+T_{36}+T_{34}+R_{16}.$$
This involution is the same as the composition $\tau\circ
p_{16}\circ p_{36}$.

\bigskip (iv) Let us consider the two pencils defined by planes
through non-skew edges. Take the edges $T_{15}$ and $T_{25}$.
Computations similar to (iii) show that the involution coincides
with the projection from $N_{56}$.

\bigskip
(v) Consider a pencil of type 3 with reducible fibres
$$C_{16}+C_{14},\quad N_{26}+N_{26}',\quad N_{45}+N_{45}',\quad
C_{23}+2T_{15}+N_{16}+N_{56}+N_{24},$$
$$2T_{36}+2T_{34}+2N_{36}+N_{13}+N_{12}+N_{46}+N_{35}.$$  We verify that
$T_{46},T_{56},T_{26}$ and $T_{25}$ are sections. Let us take
$T_{46}$ as the zero section. We have already observed that
$T_{25}$ is a 2-torsion section.  Consider the automorphism
$\phi$ of the surface defined by the inversion map with respect
to the group law on the set of sections with zero section defined
by $T_{46}$. Obviously, $T_{46}$ and $T_{25}$ are invariant.
Also, the components of reducible fibres which they intersect are
also invariant. Finally, the multiple components of fibres are
invariant. This easily shows that all irreducible components of
fibres are invariant. Let $M$ be the sublattice spanned by
irreducible components of fibres and the sections $T_{46}, T_{25}$. Its
rank is 15. Let $\alpha\in S_H$ be a primitive vector orthogonal to $M$.
Obviously, $\phi(\alpha) = -\alpha$. One easily finds the vector $\alpha$.
$$\alpha =
2\eta_H-2\eta_S+2(T_{46}+T_{25}+T_{34}+T_{36}+T_{15})+(N_{36}+N_{46}+N_{56}+N_{13}+N_{24}).$$
We leave to the reader to check that $\alpha$ is a (-6)-root of
$S_H$. Therefore, the action of $\phi$ on the Picard lattice $S_H$
is a reflection
$$\phi(x) = x-{2(x\cdot \alpha)\over \alpha^2}\alpha.\eqno (6.1)$$

\head
{\bf 7. The Leech lattice}\endhead

\bigskip\noindent
 We follow the notation and the main ideas from Kond$\bar {\roman o}$'s
paper [Ko]. First we embed the Picard lattice $S_H$ of $\tilde H$ in the
lattice $L = \Lambda\perp U \cong II^{1,25}$, where $\Lambda$ is the Leech
lattice and $U$ is the hyperbolic plane. We denote each vector
$x\in L$ by $(\lambda,m,n)$ where $\lambda\in \Lambda$, and $x =
\lambda+mf+ng,$ with $f,g$ being the standard generators of $U$.
Note that, $r = (\lambda,1,-1-{<\lambda,\lambda>\over 2})$
satisfies $r^{2} = -2$. Such vectors will be called {\it Leech roots}.
Recall that $\Lambda$ can be defined as a certain lattice in
$\bbR^{24} = \bbR^{\bbP^1(\bbF_{23})}$ equipped with inner
product $<x,y> = -{x\cdot y\over 8}$ (see [CS]). For any subset $A$ of
$\Omega = \bbP^1(\bbF_{23})$ let $\nu_A$ denote the vector
$\sum_{i\in A} e_i$, where $\{e_\infty, e_0,\ldots,e_{22}\}$ is
the standard basis in $\bbR^{24}$. Then $\Lambda$ is defined as a
lattice generated by the vectors $\nu_\Omega-4\nu_\infty$ and
$2\nu_K$, where $K$ belongs to the  SL$(2,\bbF_{23})$-orbit of
the ordered subset $(\infty,0,1,3,12,15,21,22)$ of $\Omega$.
These sets form a Steiner system $S(5,8,24)$ of eight-element
subsets of $\Omega$ such that any five-element subset belongs to
a unique element of $S(5,8,24)$. All such sets are explicitly
listed in
 [To].

\medskip
\plainproclaim Lemma 7.1. There is a primitive embedding of $S_H$ in
$L$ such that $S_H^\perp$ contains a sublattice of index $2$
isomorphic to the root lattice $R = A_5+A_1^5$.

{\sl Proof.} Consider the following Leech roots:
$$x = (4\nu_\infty+\nu_\Omega,1,2),\quad y = (4\nu_0+\nu_\Omega,
1,2),\quad z = (0,1,-1),
$$
$$x_0 = (4\nu_\infty+4\nu_0,1,1),\quad x_i = (2\nu_{K_i},1,1), \quad i =
1,\ldots,5,$$ where
$$K_1 = \{\infty,0,1,2,3,5,14,17\},$$
 $$K_2 =
\{\infty,0,1,2,4,13,16,22\}, K_3 = \{\infty,0,1,2,6,7,19,21\},$$
$$K_4 = \{\infty,0,1,2,8,11,12,18\}, K_5 = \{\infty,0,1,2,9,10,15,20\}.$$
It is easy to verify that the inner product of the vectors $x,y,z,x_i$ is
described by the following (reducible) Coxeter-Dynkin diagram:
$$\matrix x&{}&z&{}&y&{}&x_0&{}&x_1&{}&x_2&{}&x_3&{}&x_4&{}&x_5\\
\bullet&--&\bullet&--&\bullet&{}&\bullet&{}&\bullet&{}&
\bullet&{}&\bullet&{}&\bullet&{}&\bullet\endmatrix$$
Thus these vectors span a root subslattice $R_0$ of $L$ isomorphic to
$A_3+A_1^6$. We shall add one more vector to the previous set. Let
$$r_0 = (2\nu_K,1,1) \quad  K_0 = \{\infty,1,2,3,4,6,15,18\}.$$
One verifies that $<r_0,y> = <r_0,x_0> = 1$ and $<r_0,x> = <r_0,z>
= <r_0,x_i> = 0, i \ne 0$. Thus the new Dynkin diagram looks as

$$\matrix x&{}&z&{}&y&{}&r_0&{}&x_0&{}&x_1&{}&x_2&{}&x_3&{}&x_4&{}&x_5\\
\bullet&--&\bullet&--&\bullet&--&\bullet&--&\bullet&{}&\bullet&{}&
\bullet&{}&\bullet&{}&\bullet&{}&\bullet\endmatrix$$ So, the new
lattice $R$ spanned by $x,y,z,x_i,r_0$ is isomorphic to
$A_5+A_1^5$. Let
$$\theta = {1\over 2}(x+y+\sum_{i=0}^5x_i).$$
Let $q_M:D(M) = M^*/M\to \bbQ/2\bbZ$ denote  the discriminant quadratic
form of a lattice
$M$. We have
$$D(A_5) = ({5x+4z+3y+2r_0+x_0\over 6}), \quad q_{A^5} = <-{5\over 6}>,$$
$$D(A_1) = ({x_i\over 2}),\quad q_{A_1} = <-{1\over 2}>.$$ Write
$\bar\theta =
\theta\in D(R)$ as a sum of $v = {5x+4z+3y+2r+x_0\over 6}$ and $\beta_i =
{x_i\over 2}$ modulo $R$. We have $\bar\theta =
3v+\sum_{i=1}^5\beta_i$. It is checked that $\bar\theta$ spans a
subgroup $A$ of $D(R)$ which is isotropic
with respect to the discriminant quadratic form $q_R$. Let $T$ be
the overlattice of $R$ corresponding to the group $A^\perp/A$. It
is easy to check that
$$q_T=q_R|A^\perp/A=<v+\beta_1, v+\beta_2, v+\beta_3,
v+\beta_4, v+\beta_5>/\bar\theta $$ $$=
<v+\beta_1+\beta_3+\beta_5, \beta_1+\beta_2>\oplus
<\beta_1+\beta_2+\beta_3+\beta_4, \beta_1+\beta_2+\beta_4+\beta_5>
=q_{A_2(2)}\oplus q_{U(2)}.$$ Thus $T = <R,\theta>$ has the same
discriminant quadratic form as the Picard lattice $S_H$. We skip
the verification that $T$ is primitive. It is similar to the
proof of Lemma 4.1 from [Ko].
 Thus the
orthogonal complement of $T$ in $L$ is a primitive lattice of
rank 16 with the same discriminant form. Now the result follows
from the uniqueness theorem of Nikulin [Ni].

\medskip\noindent
{\bf Remark 7.2.} Note that the orthogonal complement of the
lattice $R_0$ in $L$ is isomorphic to the Picard lattice $S_K$ of
a general Jacobian Kummer surface. It contains the sublattice
$S_H$ as the orthogonal complement of the projection of $r_0$ in
$S_K$.

\medskip
One can give an explicit formula for twenty vectors
$N_\alpha',T_\beta'\in T^\perp$ whose intersection matrix
coincides with the intersection matrix of the divisor classes
$N_\alpha$, $T_\beta$. We shall identify $N_\alpha', T_\beta'$
with $N_\alpha,T_\beta$.

Explicitly, $N_\alpha,T_\beta$ correspond to the Leech vectors
$(2\nu_K,1,1)$, where
$K\subset
\Omega$ is
$$N_{45} : \{\infty,0,1,3,4,11,19,20\},\quad N_{56}:
\{\infty,0,1,3,6,8,10,13\},$$
$$N_{24} : \{\infty,0,1,3,7,9,16,18\},\quad N_{26}:
\{\infty,0,1,3,12,15,21,22\},$$
$$N_{36} : \{\infty,0,1,4,5,7,8,15\},\quad N_{35}:
\{\infty,0,1,4,6,9,12,17\},$$
$$N_{46} : \{\infty,0,1,4,10,14,18,21\},\quad N_{16}:
\{\infty,0,1,5,6,18,20,22\},$$
$$N_{13} : \{\infty,0,1,6,11,14,15,16\},\quad N_{12}:
\{\infty,0,1,13,15,17,18,19\},$$
$$T_{16} : \{\infty,0,2,3,4,8,9,21\},\quad T_{34}:
\{\infty,0,2,3,6,12,16,20\},$$
$$T_{14} : \{\infty,0,2,3,7,11,13,15\},\quad T_{36}:
\{\infty,0,2,3,10,18,19,22\},$$
$$T_{26} : \{\infty,0,2,4,5,6,10,11\},\quad T_{25}:
\{\infty,0,2,4,7,17,18,20\},$$
$$T_{15} : \{\infty,0,2,4,12,14,15,19\},\quad T_{56}:
\{\infty,0,2,5,15,16,18,21\},$$
$$T_{46} : \{\infty,0,2,6,8,15,17,22\},\quad T_{23}:
\{\infty,0,2,6,9,13,14,18\}.$$

\medskip
For Leech roots $r,r'\in L$ corresponding to the Leech vectors
$\lambda,\lambda'$ we have
$$(r,r') = \cases 0,&\text{if $\lambda-\lambda'\in \Lambda_4$};\\
1,&\text{if $\lambda-\lambda'\in \Lambda_6$}.\endcases$$
Here
$$\Lambda_4 = \{x\in \Lambda:(x,x) = 4\} =\{(\pm 2^8,0^{16}), (\pm 3,\pm
1^{23}), (\pm 4^2,0^{22})\},$$
$$\Lambda_6 = \{x\in \Lambda:(x,x) = 6\} =\{(\pm 2^{12},0^{12}), (\pm 3^3,\pm
1^{21}), (\pm 4,\pm 2^8,0^{15}), (\pm 5,\pm 1^{23})\}.$$

 Let
$$\omega = (0,0,1)\in L.$$
It is called the {\it Weyl vector} of the lattice $L$. It is an isotropic
vector with the property that
$$(\omega,l) = 1 \quad\hbox{ for any Leech root $l$}.$$

\medskip
\plainproclaim Lemma 7.3. The projection $\omega'$ of $\omega$ in
$S_H$ is equal to the vector
$$\tilde\Delta =  \calN+\calT.$$

{\sl Proof.} Note that
$(\omega',N_{\alpha})=(\omega',T_{\beta})=1$ for all $\alpha,
\beta$. On the other hand, the divisor $ \calN+\calT=\sum
N_\alpha+\sum T_\beta$ has the same property. The result follows
from Lemma 3.1.

\head
{\bf 8. Automorphisms of a general Hessian quartic surface}\endhead

Let $X$ be a K3 surface with Picard lattice $S$. The automorphism
group Aut$(X)$ of $X$ has a natural representation $\rho:\Aut(X)\to
O(S)$ in the orthogonal group of $S$. Let $W_2(S)$ denote the subgroup
generated by reflections in vectors $r$ with
$r^2 = -2$. This group is a normal subgroup of $O(X)$ and the
induced homomorphism
$$\rho:\Aut(X)\to O(S)/W_2(S)$$
has a finite kernel and a finite cokernel. This non-trivial
result follows from the Global Torelli Theorem for algebraic K3
surfaces proven by I. Shafarevich and I. Piatetski-Shapiro [PS].
Let us describe the kernel and the cokernel. First of all the
group $O(S)/W_2(S)$ has the following interpretation. Let $V_S =
\{x\in S\otimes \bbR: x^2 > 0\}$ and $V_S^+$ be its connected
component containg an ample divisor class. The group $W_2(S)$ has
a fundamental domain $P(S)$ in $V_S^+$ (a cone over a convex
polytope in the corresponding Lobachevski space $V_S^+/\bbR_+$).
It can be choosen in such a way that its bounding hyperplanes are
orthogonal to effective classes $r$ with $r^2 = -2$ and it
contains the ray spanned by an ample divisor class. Let $A(P(S))
\subset O(S)$ be the group of  symmetries of $P(S)$. Then $O(S)$
is equal to the semi-direct product $W_2(S)\rtimes A(P(S))$ of
$W_2(S)$ and $A(P(S))$. The image of $\Aut(X)$ in $O(S)$ is
contained in $A(P(S))$. Let $D(S) = S^*/S$ be the discriminant
group of the lattice and $q_S$ be the discriminant quadratic form
on $D(S)$. We have a natural homomorphism $A(P(S))\to O(q_S)$. Let
$\Gamma(S)\subset A(P(S))$ be the pre-image of $\{\pm 1\}$. Then
the image of $\Aut(X)$ in $A(P(S))$ is contained in $\Gamma(S)$
as a subgroup of index $\le 2$. It is equal to the whole group
$\Gamma(S)$ if $\Aut(X)$ contains an element acting as $-1$ on
$D(S)$. This will be the case for the Hessian surface. The kernel
of $\Aut(X)\to \Gamma(S)$ is a finite cyclic group. It is trivial
if $X$ does not admit a non-trivial automorphism preserving any
ample divisor. This happens in our case for a general Hessian
since it is known that a projective automorphism of a general
cubic surface is trivial. Thus summing up, we obtain

\medskip
\plainproclaim Proposition 8.1. Let $S_H$ be the Picard lattice of a
minimal nonsingular model $\tilde H$ of a general Hessian quartic surface
$H$. The group of automorphisms of $\tilde H$ is isomorphic to
the group
$\Gamma(S_H)$ of symmetries of $P(S_H)$ which act as $\pm 1$ on the
discriminant group of $S_H$.

\medskip
 Let $W_{\Lee}(L)$ be the subgroup of $O(L)$ generated
by reflections in Leech roots. Let $P(L)$ be its fundamental
domain in the $V_L^+$ where $V_L^+$ is one of the components of
$V = \{x\in L\otimes \bbR:x^2 > 0\}$ which contains the Weyl
vector $\omega$. By a result of J. Conway, $O(L)$ is equal to the
semi-direct product $W_{\Lee}(L)\rtimes A(P(L))$ of $W_{\Lee}(L)$ and the
group of symmetries $A(P(L))$ of $P(L)$, and the latter is isomorphic to
the group $\Lambda\rtimes O(\Lambda)$ of affine automorphisms of
$\Lambda$.

\medskip
Now put
$$P(S_H)' = P(L)\cap V(S_H)^+.$$

It is known that $P(S_H)'$ is non-empty, has only finitely many
faces and of finite volume (in the Lobachevski space). Also it is
known ([Bo], Lemmas 4.1-4.3) that the projection $\omega'$ of the
Weyl vector is contained in $P(S_H)'$. Applying Lemma 7.3, we see
that $P(S_H)'$ contains an ample divisor class and hence
$P(S_H)'$ is a part of $P(S_H)$.

\plainproclaim Lemma 8.2. Let $G = \Aut(P(S_H)')\subset
O(S_H)$ be the group of symmetries of $P(S_H)'$. Then $G$ is
isomorphic to $\bbZ/2\times S_5$ and can be realized as the
subgroup of $O(S_H)$ generated by the Enriques involution and the
group of symmetries of the Sylvester pentahedron. It extends to a
subgroup of $O(L)$ which leaves the root lattice $R$ invariant
and induces the isometries of $R$ defined by the symmetries of
the Dynkin diagram of $R$.

{\sl Proof.} This is almost word by word repetition of the proof of
Lemma 4.5
in [Ko].

\medskip
To find a certain set of generators of $\Aut(\tilde H)$
containing the Enriques involution, we use the following strategy
suggested in [Ko]. First we shall find the hyperplanes which
bound $P(S_H)'$. They correspond to rank 11 root sublattices $R'$
generated by $R$ and some Leech root $r$. The hyperplane $\{x\in
V_L^+: (x,r) = 0\}$ is a boundary wall of $P(L)$ and has non-empty
intersection with $P(S_H)$. Then for such a hyperplane we find an
automorphism of $\tilde H$ which maps one of the two half-spaces
defined by this hyperplane to the opposite one. Let $N$ be the
group generated by these automorphisms. Then we check that for any
automorphism $\gamma\in \Aut(\tilde H)$ one can find $\delta\in N$
such that $\delta\circ \gamma$ is a symmetry of $P(S_H)'$.
Applying Lemma 8.2, we conclude that $\delta\circ \gamma$ is
either identity or the Enriques involution (since the latter and
the identity are the only elements of $\Aut(P(S_H)')$ which act as
$\pm 1$ on the discriminant group of $S_H$).

\medskip
The next lemma  is a simple repetition of the computations
from [Ko], Lemma 4.6.

\bigskip \plainproclaim Lemma 8.3. Let $r$ be a Leech root. Assume
that $r$ and $R$ generate a root lattice $R'$ of rank 11. Then
one of the following cases occurs:
\item{(0)} $R' = A_5\oplus A_1^6$, where $r$ is orthogonal to $R$;
\item{(1a)} $R' = D_6\oplus A_1^5$, where $ (r,r_0) = 0$, $(r,z) =
1$;
\item{(1b)} $R' = D_6\oplus A_1^5$, where $ (r,r_0) = 1$, $(r,z) =
0$;
\item{(2)} $R' = A_1^3\oplus A_3\oplus A_5$, where $(r,x_i) = (r,x_j) = 1, i,j
\ne 0$;
\item{(3a)} $R' = A_7\oplus A_1^4$, where $(r,x_0) = (r,x_i) = 1$;
\item{(3b)} $R' = A_7\oplus A_1^4$, where $(r,x) = (r,x_i) = 1$.

{\it Moreover, in case (0), $r$ is one of the twenty Leech roots
corresponding to $N_\alpha$ and $T_\beta$}.

{\it In case (1a) up to a symmetry of $P(S_H)'$ we can choose
$r=(\lambda,1,-1-{<\lambda,\lambda>\over 2})$  corresponding to
the Leech vector
$$\lambda=(\xi_\infty,\xi_0,\xi_{j_1},\xi_{j_2},\ldots,\xi_{j_6},\xi_{j_7},\ldots,\xi_{j_{22}}) =
(3,3,3,-1,\ldots,-1,1,\ldots,1),$$ where $K =
\{\infty,0,j_1,\ldots,j_6\}$ is an octad satisfying $|K\cap K_0| =
\{\infty, j_1\}$, $|K\cap K_i| = 4$ and $j_1\in K_i$ for
$i=1,\ldots,5$}.

{\it In case (1b) $r$ corresponds to $N_\alpha, \alpha\in
\{0,14,15,23,25,34\}$ or to $T_\beta, \beta\in
\{0,12,$
$13,24,35,45\}.$ Together with ``old '' $N_\alpha$ and
$T_\beta$ they define 32 vectors spanning the Kummer overlatice
$S_K$ of $S_H$. }

{\it In case (2) $r$ corresponds to a Leech vector
$$\lambda=2\nu_K, \quad 0,\infty\in K, \quad |K\cap K_i| = |K\cap K_j| = 2, $$
$$|K\cap
K_l| = 4\quad
\hbox{ for $l=0,\ldots,5, l\ne i,j$}.$$}

{\it In case (3a), if $r$ meets $x_0,x_i$, it corresponds to a
Leech vector
$$\lambda=\nu_\Omega-4\nu_k, \quad k\ne 0,\infty, k\in K_i, k\not\in K_j, j=0,1,\ldots,5, j\ne i.$$}

{\it In case (3b), if $r$ meets $x,x_i$, it corresponds to a
Leech vector
$$\lambda=4\nu_0-\nu_K+\nu_{\Omega-K}, \quad  0\in K, \infty\not\in K,$$
$$|K\cap K_0| =0, |K\cap K_i| = 4, |K\cap K_l| = 2, l=1,\ldots,5,
l\ne i.$$}

\medskip\noindent
{\bf Remark 8.4}. The number of vectors $r$ in
case (1a) is equal to $12$. They correspond to the following
octads $K$:
$$\{\infty,0,1,5,9,11,13,21\},\quad \{\infty,0,1,7,10,11,17,22\},\quad
\{\infty,0,1,7,12,13,14,20\},$$
$$\{\infty,0,1,8,9,14,19,22\},\quad \{\infty,0,1,8,16,17,20,21\},\quad
\{\infty,0,1,5,10,12,16,19\},$$
$$\{\infty,0,2,5,7,9,12,22\},\quad \{\infty,0,2,5,8,13,19,20\},\quad
\{\infty,0,2,7,8,10,14,16\},$$
$$\{\infty,0,2,9,11,16,17,19\},\quad \{\infty,0,2,10,12,13,17,21\},\quad
\{\infty,0,2,11,14,20,21,22\}.$$

\noindent The number of vectors $r$ in case (2) is equal to $10 =
{5\choose 2}$. When $i = 1, j = 2$, $r$ corresponds to the
following octad $K$:
$$\{0,\infty,6,7,10,12,15,18\}.$$

\noindent The number of vectors $r$ in case (3a) is equal to
$5\cdot3 = 15$. When $i = 1$, they correspond to the following
Leech vectors:
$$\nu_\Omega-4\nu_k,\quad  k\in \{5,14,17\}.$$

\noindent The number of vectors $r$ in case (3b) is equal to
$5\cdot3 = 15$. When $i = 1$, they correspond to the following
octads $K$:
$$\{0,5,9,12,13,14,17,19\},\quad \{0,5,10,11,14,16,17,21\},\quad
\{0,5,7,8,14,17,20,22\}.$$

\bigskip\plainproclaim Lemma 8.5. Let $r$ be a Leech root as in the
previous lemma. Let $r = r_1+r_2$, where $r_1\in S_H^*$ and
$r_2\in T^*$. Then

{\it Case $(0):$ $r_1 = r$};

{\it Case $(1a):$ $r_1 = r+{1\over 3}(2x+4z+3y+2r_0+x_0)$,
$(r_1,r_1) = -2/3$;}

{\it Case $(1b):$ $r_1 = r+{1\over 3}(x+2z+3y+4r_0+2x_0)$,
$(r_1,r_1) = -2/3$;}

{\it Case $(2):$ $r_1 = r+{1\over 2}(x_i+x_j)$, $(r_1,r_1) = -1$;}

{\it Case $(3a):$ $r_1 = r+{1\over 6}(x+2z+3y+4r_0+5x_0)+{1\over
2}x_i$, $(r_1,r_1) = -2/3$;}

{\it Case $(3b):$ $r_1 = r+{1\over 6}(5x+4z+3y+2r_0+x_0)+{1\over
2}x_i$, $(r_1,r_1) = -2/3$}.

\bigskip
{\sl Proof.} Case (0) is obvious.

Case (1a): Let $x^*,y^*,z^*,r_0^*,x_i^*$ denote the dual basis of
$x,y,z,r_0,x_0,\ldots,x_5$. Since $(r,z) = 1, (r,x)= (r,y) =
(r,r_0) = (r,x_i) = 0$, we see that $r_2 = z^* =
-(2x+4z+3y+2r_0+x_0)/3$ and hence $(r_1,r_1) = (r-r_2,r-r_2) =
(r,r)-(r_2,r_2) = -2+4/3 = -2/3.$

Case (2): Here $r_2 = x_i^*+x_j^* = -(x_i+x_j)/2$. This gives
$(r_1,r_1) =-1.$

Case (3a): $r_2 = x_0^*+x_i^* = -{1\over
6}(x+2z+3y+4r_0+5x_0)-{1\over 2}x_i$. This gives $(r_1,r_1)=-2/3$.

\bigskip
Let $s_{r}: v\to v+(v,r)r$ be the reflection of $L$ in a Leech
root $r$. Then the restriction of $s_r$ to $S_H\otimes \bbQ$ is a
reflection
$$s_{r_1}(v) =  v+(v,r_1)r_1,$$
where $r_1$ is the projection of $r$ onto $S_H\otimes \bbQ$. This
is, in general, not an isometry of $S_H$.

To find an automorphism of $\tilde H$ corresponding to each
vector in case (1a)-(3b), we need to express $r_1$ in terms of
$N_\alpha, T_\beta$.

\bigskip
{\bf Case (1a).} Let  $r$ correspond to the Leech vector
$$3\nu_0+3\nu_\infty+3\nu_1-\nu_5-\nu_{9}-\nu_{11}-\nu_{13}-\nu_{21}+\nu_{\Omega-K},$$
where $K = \{0,1,\infty,5,9,11,13,21\}.$
It is immediately verified that
$$(r, N_\alpha) = 0, \quad\hbox{for all $\alpha$},$$
$$(r, T_\beta) = \cases
1,&\text{if $\beta = 16,26,56,14,23$};\\
 0,&\text{otherwise}.\endcases$$
This
determines
$r_1$ in the form of a linear combination of $N_\alpha, T_\beta$.
Writing down the corresponding system of linear equations and
solving it we obtain
$$r_1 ={1\over
15}(-2\calT+2\calN+10(T_{36}+T_{46}+T_{15}+T_{25}+
T_{34})+5(N_{36}+N_{46}+N_{56}+N_{13}+N_{24})).\eqno (8.1)$$
 Let
$$\alpha = 3r_1.$$
By (3.1), $\calN-\calT = 5\eta_H-5\eta_S$, so that
$$\alpha = 2\eta_H-2\eta_S+2(T_{36}+T_{46}+T_{15}+T_{25}+T_{34})+(N_{36}+N_{46}+N_{56}+N_{13}+N_{24}).$$
It is easy to check that  $\alpha$ is a primitive vector of $S_H$
and a (-6)-root.

\bigskip
{\bf Case (1b).} These are conjugate to those in case (1a) by the
Enriques involution $\tau$ which can be viewed as an automorphism
of $P(S_H)'$. For example, take $r = T_{12}$. We know that
$(r_1,T_\alpha) = 0$ for all $\alpha$ and
$$(r_1,N_\beta) = \cases 1,&\text{if $\beta = 16,26,12,35,45$};\\
0,&\text{otherwise}.\endcases$$
Applying the Enriques involution $\tau$,
we see that $(\tau(r_1),N_\alpha) = 0$ for all $\alpha$ and
$$(\tau(r_1),T_\beta) =\cases 1,&\text{if $\beta = 16,26,56,14,23$};\\
0,&\text{otherwise}.\endcases$$ Thus, $\tau(r_1)$ is nothing but
the vector appeared in Case (1a) above.

\bigskip
{\bf Case (2).} Let $r = (2\nu_K,1,1)$, where $K =
\{0,\infty,6,7,10,12,15,18\}$. We have
$$(r,N_\alpha) =\cases 1,&\text{if $\alpha = 45$};\\
0,&\text{otherwise},\endcases$$
$$(r,T_\beta) =\cases 1,&\text{if $\beta = 16$};\\
0,&\text{otherwise}.\endcases$$ Notice that $N_{45}$ is the
vertex opposite to the edge $T_{16}$. It is easy to see that
$(r,N_\alpha) = ({1\over 2}(C_{16}-N_{45}),N_\alpha)$ for all
$\alpha$ and $(r,T_\beta) = ({1\over 2}(C_{16}-N_{45}),T_\beta)$
for all $\beta$. Thus
$$r_1 = {1\over 2}(C_{16}-N_{45}) = {1\over
2}(\eta_H-2T_{16}-N_{16}-N_{12}-N_{13}-N_{45}).\eqno (8.2)$$ Note that $\alpha
= 2r_1$ is a $(-4)$-root of $S_H$. Also note that $\tau(\alpha)=\alpha$.

\bigskip
{\bf Case (3a).} Let $r$ correspond to the Leech vector
$\nu_\Omega-4\nu_5$. We have
$$(r,N_\alpha) =\cases 1,&\text{if $\alpha = 16,36$};\\
0,&\text{otherwise},\endcases$$
$$(r,T_\beta) =\cases 1,&\text{if $\beta = 26,56$};\\
0,&\text{otherwise}.\endcases$$
By solving a system of linear equations we
find that
$$r_1={1\over 2}(N_{26}+N_{45}+N_{56}+N_{24})-{1\over
3}(N_{16}+N_{36})+$$
$${1\over 3}(N_{13}+N_{46}) +{2\over
3}(T_{46}+T_{25})+{1\over 3}(T_{15}+T_{23}).\eqno (8.3)$$ It can be
checked that the minimum positive integer $k$ with $kr_1\in S_H$ is 6, so
that
$$6r_1=3(N_{26}+N_{45}+N_{56}+N_{24})-2(N_{16}+N_{36})+2(N_{13}+N_{46})$$
$$+4(T_{46}+T_{25})+2(T_{15}+T_{23})$$ is a primitive (-24)-vector in $S_H$. We
remark that $\alpha=6r_1$ is NOT a root.

\bigskip
{\bf Case (3b).} These are conjugate to those in case (3a) by the
Enriques involution $\tau$. For example, let $r$ correspond to the
Leech vector
$$4\nu_0-\nu_K+\nu_{\Omega-K},$$
where $K=\{0,5,9,12,13,14,17,19\}$. We have
$$(r,N_\alpha) =\cases 1&\text{if $\alpha = 12,35$};\\
0&\text{otherwise.}\endcases$$
$$(r,T_\beta) =\cases 1,&\text{if $\beta = 15,23$};\\
0&\text{otherwise.}\endcases$$
Applying $\tau$, we see that
$$(\tau(r_1),T_\alpha) =\cases 1,&\text{if $\alpha = 26,56$};\\
0,&\text{otherwise},\endcases$$
$$(\tau(r_1),N_\beta) =\cases 1,&\text{if $\beta = 16,36$};\\
0&\text{otherwise.}\endcases$$ and that $\tau(r_1)$ is the one appeared
in Case (3a) above.

\bigskip
Let us find an involution $\sigma$ of $\tilde H$ corresponding to
a hyperplane defined by a Leech root $r$ of type $(1a), (2)$, and
$(3a)$. If $r_1$ denotes the projection of $r$ to $S_H$ we need
that
$$\sigma(r_1) = -r_1.$$

\bigskip {\bf Automorphisms of type (1)}
\smallskip {\bf (1a)} This is the
inversion of an elliptic pencil of type 3. For example, the
involution $\phi$ from (6.1) corresponds to the vector
$r_1$ from (8.1). Since $\alpha = 3r_1$ is a (-6)-root, we get
$\phi(r_1)=-r_1$. Also
$$\phi(\omega') = \omega'+3(\omega',r_1)r_1 = \omega'+15r_1,\eqno (8.4)$$
where $\omega'=\sum_\alpha N_\alpha +\sum_\beta T_\beta$ is the
projection of the Weyl vector(see Lemma 7.3).

Let us explain why the number of such
automorphisms is 12. The inversion of the same elliptic fibration
with respect to the different zero section $T_{26}$ (and a
2-torsion $T_{56}$) gives another involution. Since there are 30
elliptic fibrations of type 3, we get in this way 60 involutions but only
12 of them are different.

\smallskip {\bf (1b)} This is the conjugate involution
$\tau\circ\phi\circ\tau$.

\bigskip
{\bf Automorphisms of type (2)} \medskip These are the 10
projections $p_{\alpha}$ from a node $N_{\alpha}$. For example,
consider $p_{45}$ whose action can be computed as in Section
6-(ii), and take $r_1$ computed in (8.2). It is easy to check that
 $p_{45}(r_1)=-r_1$. Note that $p_{45}$ acts nontrivially on the
hyperplane $r_1^{\perp}$ in
$V(S_H)^+$. In other words, $p_{45}$ is not a reflection but works like a
reflection. Also observe  that
$$p_{45}(\omega')=\omega'+4r_1=\omega'+2(\omega',r_1)r_1.\eqno (8.5)$$

\bigskip {\bf Automorphisms of type (3)}

\smallskip {\bf (3a)} This is the inversion of an elliptic pencil of type
2. To see this, take $r_1$ computed in (8.3).

 Consider the elliptic fibration
$$|F_{N_{12},T_{26}}|=|F_{N_{35},T_{56}}|=|T_{23}+N_{26}+T_{25}+N_{56}+T_{15}+N_{24}+T_{46}+N_{45}|
= $$
$$|C_{56}+T_{34}+N_{12}+T_{16}|=|C_{26}+T_{36}+N_{35}+T_{14}|.$$
Observe that $N_{13}$, $N_{46}$, $N_{16}$, $N_{36}$ are sections.
Take $N_{13}$ as the zero section. Then $N_{46}$ is a 2-torsion.
Let $f_r=f_{N_{16},N_{36},T_{26},T_{56}}$ be the inversion. Take
$$D_1=(N_{26}+N_{45})+2(N_{56}+N_{24})-N_{36}+(N_{13}+N_{46})
+2(T_{46}+T_{25}+T_{15}),$$
$$D_3=2(N_{26}+N_{45})+(N_{56}+N_{24})-N_{16}+(N_{13}+N_{46})
+2(T_{46}+T_{25}+T_{23}).$$ Then $D_1$ and $D_3$ are effective
(-2)-vectors, and
$$D_1+N_{16}-2N_{13}\sim
T_{36}-T_{14}+T_{23}+2N_{26}+3T_{25}+2N_{56}+T_{15}+N_{24}
+T_{46}+N_{45},$$
$$D_3+N_{36}-2N_{13}\sim
T_{16}-T_{34}+T_{23}+2N_{26}+3T_{25}+2N_{56}+T_{15}+N_{24}
+T_{46}+N_{45}.$$
We see that $D_1$ and $D_3$ are sections and
$D_1\oplus N_{16}=D_3\oplus N_{36}=0$. These determine the action
of $f_r$ as follows:
$$N_{16}\longleftrightarrow D_1, \quad N_{26}\longleftrightarrow N_{56},\quad
N_{36}\longleftrightarrow D_3,\quad N_{46}\longleftrightarrow N_{46},$$
$$N_{12}\longleftrightarrow C_{56}, \quad N_{13}\longleftrightarrow N_{13},
\quad N_{24}\longleftrightarrow N_{45}, \quad
N_{35}\longleftrightarrow C_{26},$$
$$T_{15}\longleftrightarrow T_{23},$$
$$f_r(T_{\beta})=T_{\beta}, \quad \beta=16, 36, 46, 14, 25, 34,$$
$$f_r(T_{26})= 4\eta_H
-2(T_{26}+T_{56})-T_{26}-T_{15}-T_{23}-2(N_{16}+$$
$$N_{26}+N_{36}+N_{56}+N_{12}+N_{24}+N_{35}+N_{45}),$$
$$f_r(T_{56})= 4\eta_H
-2(T_{26}+T_{56})-T_{56}-T_{15}-T_{23}-2(N_{16}+$$
$$N_{26}+N_{36}+N_{56}+N_{12}+N_{24}+N_{35}+N_{45}).$$
Note that $$D_1+D_3=N_{16}+N_{36}+6r_1,$$
so that
$f_r(6r_1)=-6r_1.$ Again
$f_r$ is not a reflection, but works like a reflection.

\med  The involution $f_r$ is not symmetric in the sense that
$$f_r(T_{26}+T_{56})\neq T_{26}+T_{56}+6r_1.$$
On the other hand, the map
$g_r=g_{N_{16},N_{36},T_{26},T_{56}}:=f_{N_{16},N_{36},T_{26},T_{56}}\circ
p_{35}\circ p_{12}$ is symmetric. We have

$$\left[\matrix N_{16}& N_{26}& N_{36}& N_{46}& N_{56}& N_{12}& N_{13}& N_{24}&
N_{35}& N_{45}\\D_3& N_{45}& D_1& N_{13}& N_{24}& N_{12}& N_{46}&
N_{56}& N_{35}& N_{26}\endmatrix\right],$$

\med
$$\left[\matrix T_{16}& T_{26}& T_{36}& T_{46}& T_{56}& T_{14}& T_{15}& T_{23}&
T_{25}& T_{34}\\T_{34}& G_2& T_{14}& T_{25}& G_5& T_{36}& T_{15}&
T_{23}& T_{46}& T_{16}\endmatrix\right],$$ where $g_r$ interchanges
two elements in the same column, and
$$G_2=T_{26}+T_{15}+T_{23}+2(T_{46}+T_{25})+2(N_{26}+N_{24})-(N_{16}+N_{36})+N_{56}+N_{45}+N_{13}+N_{46},$$
$$G_5=T_{56}+T_{15}+T_{23}+2(T_{46}+T_{25})+2(N_{56}+N_{45})-(N_{16}+N_{36})+N_{26}+N_{24}+N_{13}+N_{46}.$$

In fact, this is an isometry of $S_H$ acting as $-1$ on the
discriminant form of $S_H$, and hence sending $6r_1$ to $-6r_1$.
In fact, following an idea from [Ke2], we found first the lattice
involution $g_r$ and then $f_r$ which relaizes it geometrically.

\med  Since $g_r(N_{16}+N_{36})=N_{16}+N_{36}+6r_1$ and
$g_r(T_{26}+T_{56})= T_{26}+T_{56}+6r_1$, we have
$$g_r(r_1) =
-r_1, \quad g_r(\omega')=\omega'+12r_1=\omega'+3(\omega',r_1)r_1,\eqno
(8.6)$$ where $\omega'=\sum_\alpha N_\alpha +\sum_\beta T_\beta$ is the
projection of the Weyl vector.

\med
{\bf (3b)} This is the involution $\tau\circ g_{\tau(r)}\circ \tau.$

\med {\bf Remark 8.6.}  The translation by the 2-torsion $N_{46}$
is the same as $\tau\circ$(involution defined by the skew lines
$T_{26}, T_{56}$) = $p_{35}\circ p_{12}$.

\med {\bf Remark 8.7.} If we take $N_{16}$ as the zero section
(and $N_{36}$ a 2-torsion), then the inversion map corresponds to
a Leech root of type (3b), more precisely, to $r$ with
$$(r,N_\alpha) =\cases 1,&\text{if $\alpha = 46,13$};\\
0,&\text{otherwise.}\endcases$$
$$(r,T_\beta) =\cases 1,&\text{if
$\beta = 26,56$};\\ 0,&\text{otherwise.}\endcases$$

\medskip \plainproclaim Theorem 8.8. The automorphism group of $\tilde H$
is generated by the Enriques involution $\tau$, the $10$
projections $p_{\alpha}$, the $15$ inversion automorphisms $f_r$
of elliptic pencils of type 2, and $12$ inversion automorphisms
$\phi_r$ of elliptic pencils of type 3.

{\sl Proof.} Let $\sigma_r$ be the involution corresponding to a
Leech root $r$ of type (1a), (1b), (2), (3a), or (3b), i.e.
$\sigma_r$ is one of the twelve $\phi_r$, twelve
$\tau\circ\phi_{\tau(r)}\circ\tau$, ten $p_r$, fifteen $g_r$, or
fifteen $\tau\circ g_{\tau(r)}\circ\tau$. Let $N$ be the subgroup
of Aut$(\tilde H)$ generated by them.

As we explained before (after Lemma 8.2), the result follows from the
following:

\bigskip \plainproclaim Lemma 8.9. Let $\gamma$ be an isometry of
the Picard lattice $S_H$ which preserves $P(S_H)$. Then there
exists an element $\delta\in N$ such that
$\delta\circ\gamma\in\Aut(P(S_H)')$.

{\sl Proof.} This is similar to the proof of Lemma 7.3 in [Ko].
Take an element $\delta\in N$ which realizes
$\min \{(\delta(\gamma(\omega')), \omega') : \delta\in N\}$. Then for any
$r$,
$$(\delta\circ\gamma(\omega'), \omega')\le (\sigma_r\circ\delta\circ\gamma(\omega'),
\omega')=(\delta\circ\gamma(\omega'), \sigma_r(\omega')).$$

If $\sigma_r=\phi_r$ (of type (1a)), then, applying (8.4), we get
$$(\delta\circ\gamma(\omega'), \omega')\le (\delta\circ\gamma(\omega'),
\omega')+15(\delta\circ\gamma(\omega'), r_1).$$ This means that
$(\delta\circ\gamma(\omega'), r_1)\ge 0$. Since $\omega'$ is an
interior point of $P(S_H)'$, the last inequality is strict.

If $\sigma_r=\tau\circ\phi_{\tau(r)}\circ\tau$ (of type (1b)),
then, since $\tau$ preserves $\omega'$, we have
$$\sigma_r(\omega')=\tau(\omega'+15\tau(r)_1)=\omega'+15r_1$$
and hence again $(\delta\circ\gamma(\omega'), r_1)> 0$.

The remaining cases can also be easily handled in this way by
using (8.5), (8.6). Now $(\delta\circ\gamma(\omega'), r_1)> 0$
for all $r$, so $\delta\circ\gamma(\omega')\in P(S_H)'$.

\bigskip\noindent
 {\bf Remark 8.10.} As a general Hessian quartic surface
degenerates to a general Jacobian Kummer surface, the involutions
of type (1) are refined to become projections and correlations of
the Jacobian Kummer, the involutions of type (2) become Cremona
transformations related to G\"opel tetrads; on Jacobian Kummer
surface there are 10 G\"opel tetrads having exactly 3 elements in
common with the hexad (the projection of $r_0$ onto the Picard
lattice of the Jacobian Kummer surface is ${1\over
4}(3\eta_{K}-2\roman{hexad})$ and the projection of a Leech root
of type (2) is ${1\over 2}(\eta_{K}-\text{G\"opel tetrad})$, and
these two vectors are orthogonal to each other), and the
involutions of type (3a) also become Cremona transformations
related to G\"opel tetrads; there are 15 G\"opel tetrads having
exactly 2 elements in common with the hexad. It looks quite
complicated to recover the involution $f_r$ or $g_r$ from its
corresponding Cremona transformation on the Jacobian Kummer
surface. Finally, the 192 new automorphisms  [Ke2] on the Jacobian
Kummer surface do not correspond to any automorphisms of the
Hessian, in other words, as a general Hessian degenerates  to a
Jacobian Kummer surface, generators of type (3b) are replaced by
new generators, which are the 192 automorphisms. Finally note that in our case all the generators send one of the half-spaces defined by the corresponding hyperplane to the opposite half-space. The new automorphisms of the Jacobian Kummer do not act in this way; they send the half-space to a half space defined by the hyperplane corresponding to the inverse automorphism.

\vglue 1in
I. Dolgachev:
 Department of Mathematics,
 University of Michigan,
Ann Arbor, MI 48109;
email: idolga{\@}umich.edu

J. Keum: Korea Institute for Advanced Study, 207-43 Cheongryangri-dong,
Dongdaemun-gu, Seoul 130-012, Korea;
email:
 jhkeum{\@}kias.re.kr

\end{document}